\documentclass[12pt,a4paper]{aclart}

\usepackage{mathptmx}
\title{Mesures et \'equidistribution sur les espaces de Berkovich}
\author{Antoine Chambert-Loir}
\address{IRMAR, Universit\'e de Rennes~1 \\
Campus de Beaulieu \\ 35042 Rennes Cedex}
\email{antoine.chambert-loir@univ-rennes1.fr}
\begin{abstract}
Les d\'emonstrations par Ullmo et Zhang de la conjecture
de Bogomolov sur les points de petite hauteur dans les vari\'et\'es
ab\'eliennes utilisent de fa\c{c}on cruciale une propri\'et\'e d'\'equidistribution
des {\og petits points\fg} dans la vari\'et\'e ab\'elienne complexe
associ\'ee. Cet article \'etudie des propri\'et\'es d'\'equidistribution
analogues aux places finies. Nos r\'esultats
s'expriment naturellement dans le cadre des espaces 
analytiques d\'efinis par Berkovich.
Le premier, en dimension quelconque,
est restreint aux {\og m\'etriques alg\'ebriques\fg};
le second n'est valable que pour les courbes mais s'applique
\`a des m\'etriques plus g\'en\'erales, notamment aux hauteurs
normalis\'ees par un syst\`eme dynamique.
\end{abstract}
\begin{altabstract}
The proof by Ullmo and Zhang of Bogomolov's conjecture about points of
small height in abelian varieties made a crucial use of an equidistribution
property for ``small points'' in the associated complex abelian variety.
We study the analogous equidistribution property at $p$-adic places. 
Our results can be conveniently stated within the framework
of the analytic spaces defined by Berkovich.
The first one is valid in any dimension
but is restricted to ``algebraic metrics'', the
second one is valid for curves, but allows for more general metrics,
in particular to the normalized heights 
with respect to dynamical systems.
\end{altabstract}

\def\Div{{\operatorname{Div}}}
\def\gm{{\mathbf G_{\mathrm m}}}
\def\C{{\mathbf C}}
\def\P{{\mathbf P}}
\def\Q{{\mathbf Q}}
\def\N{{\mathbf N}}
\def\R{{\mathbf R}}
\def\Z{{\mathbf Z}}
\def\sp{\operatorname{sp}}
\let\ra\rightarrow
\let\eps\varepsilon
\let\phi\varphi
\let\bar\overline

\def\resp{\emph{resp.}\xspace}
\def\loccit{\emph{loc.cit.}\xspace}
\def\cf{\emph{cf.}\xspace}
\def\an{{\text{\upshape an}}}
\def\abs#1{\left|{#1}\right|}
\def\norm#1{\left\| {#1} \right\|}

\def\hc{{\widehat c}}

\def\hdiv{\mathop{\widehat{\operatorname{div}}}}

\def\Gal{\operatorname{Gal}}
\def\Spf{\operatorname{Spf}}
\def\Spec{\operatorname{Spec}}

\def\div{\operatorname{div}}
\def\Pic{\operatorname{Pic}}
\def\bPic{\overline{\Pic}}
\def\alg{{\text{\upshape alg}}}
\def\bPicalg{\overline{\Pic}_{\text{\upshape alg}}}
\def\bPicfor{\overline{\Pic}_{\text{\upshape for}}}
\def\bPicint{\overline{\Pic}_{\text{\upshape int}}}
\def\tube#1{\mathopen]#1\mathclose[}
\def\sozat{\,;\,}

\begin{document}
\maketitle
\tableofcontents

\section{Introduction}

Soit $F$ un corps de nombres et soit $X$ un $F$-sch\'ema
projectif lisse.
Soit $v$ une place de~$F$ et soit $K$ le compl\'et\'e de~$F$ en~$v$,
soit $\bar K$ le compl\'et\'e d'une cl\^oture alg\'ebrique de~$K$. 
On supposera fix\'e un plongement d'une cl\^oture alg\'ebrique de~$F$ 
dans~$\bar K$.

Si~$v$ est archim\'edienne, on notera $X_v$ la vari\'et\'e analytique complexe $X(\bar K)$.
Si $v$ est ultram\'etrique, on notera $X_v$
l'espace de Berkovich 
de~$X\otimes K$, \cf~\cite{berkovich1990}; 
contentons-nous pour l'instant de dire qu'il s'agit d'une notion raisonnable 
d'espace analytique sur le corps~$K$ attach\'ee \`a~$X\otimes K$;
en particulier, $X_v$ est un espace compact, m\'etrisable,
et est connexe par arcs, si $X$ est g\'eom\'etriquement connexe.

Tout point~$x$ de $X(K)$ d\'efinit un point de $X_v$, que l'on notera
par la m\^eme lettre; on notera $\mu_x$ la mesure (de Dirac) sur~$X_v$
de masse~$1$ support\'ee en~$x$.
Plus g\'en\'eralement, si $\{x_1,\dots,x_d\}$
est l'orbite sous $\Gal(\bar K/K)$ d'un point $x \in X(\bar K)$,
on notera $\mu_x$ la mesure de probabilit\'e sur~$X_v$, somme des masses
de Dirac en les~$x_i$ divis\'ee par~$d$. (\emph{Stricto sensu,}
si $v$ est ultram\'etrique,
cette mesure est d\'efinie d'abord sur un $K'$-espace de Berkovich $X_{v'}$,
associ\'e \`a une extension finie~$K'$ de~$K$;
on consid\`ere son image directe sous le morphisme
canonique $X_{v'}\ra X_v$.)

La question de l'\'equidistribution
d'une suite de points $(x_n)$ de $X(\bar F)$
consiste \`a \'elucider le comportement des mesures $\mu_{x_n}$
lorsque $n$ tend vers l'infini.
On supposera toujours que cette suite est \emph{g\'en\'erique,}
c'est-\`a-dire  qu'une sous-vari\'et\'e donn\'ee $Z\subsetneq X$
ne contient qu'un nombre fini de termes de la suite.

Dans le contexte de la th\'eorie d'Arakelov, on s'int\'eresse
notamment aux suites de petits points.
Soit $\bar L$ un fibr\'e en droites sur~$X$ muni d'une m\'etrique
ad\'elique int\'egrable au sens de Zhang~\cite{zhang95b} ;
rappelons que cela signifie qu'en dehors d'un nombre fini de places,
cette m\'etrique est d\'efinie par un mod\`ele, et qu'en les places
restantes, elle est quotient de deux limites uniformes de m\'etriques
amples (en particulier, lisses aux places archim\'ediennes).
Comme le montre Zhang dans l'article cit\'e,
ces m\'etriques donnent lieu \`a un calcul de nombres d'intersection
arithm\'etique, d\'efini par passage \`a la limite \`a partir
de la th\'eorie d'Arakelov usuelle (\cf par exemple~\cite{gillet-s90}
ou, en ce qui concerne plus sp\'ecifiquement l'application
de la th\'eorie d'Arakelov aux hauteurs, \cite{bost-g-s94}).
En particulier, si $d=\dim X$, la formule
\[ h_{\bar L}(X) = \frac{(\hc_1(\bar L)^{d+1}|X)}{(1+d)c_1(L)^d}, \]
d\'efinit ce qu'on appelle la \emph{hauteur de~$X$} relativement
au fibr\'e m\'etris\'e~$\bar L$,
en tout cas si $c_1(L)^d>0$.

Les r\'esultats  de Szpiro-Ullmo-Zhang dans~\cite{szpiro-u-z97}
am\`enent alors la question, dans laquelle il vaut mieux supposer
que $L$ est ample.
\emph{Si $(x_n)$ est une suite g\'en\'erique de points de $X(\bar F)$
tels que $h_{\bar L}(x_n)\ra h_{\bar L}(X)$, est-il vrai
que pour toute place~$v$,
les mesures $\mu_{x_n}$ convergent faiblement ?}
Le but de cette partie est de donner des exemples significatifs
o\`u cette question a une r\'eponse positive,
voir les th\'eor\`emes~\ref{theo.berk-ample} et~\ref{theo.courbes}.
Plus g\'en\'eralement, notre approche permet d'aborder
la question analogue de l'\'equidistribution de mesures de la forme~$\mu_x$
dans un produit (fini) d'espaces~$X_v$.

L'approche que nous suivons 
est h\'erit\'ee de l'article~\cite{szpiro-u-z97}
de Szpiro, Ullmo et Zhang. Elle fait usage d'une in\'egalit\'e
fondamentale, cons\'equence du th\'eor\`eme de Hilbert-Samuel arithm\'etique:
sous certaines conditions sur le fibr\'e en droites m\'etris\'e~$\bar L$, 
on a
\[ \liminf h_{\bar L}(x_n) \geq h_{\bar L}(X) \]
pour toute suite g\'en\'erique $(x_n)$ de points de $X(\bar F)$.
Le th\'eor\`eme de Hilbert-Samuel arithm\'etique a d'abord \'et\'e d\'emontr\'e
par Gillet-Soul\'e~\cite{gillet-s88,gillet-s92}
dans le cadre de la th\'eorie d'Arakelov {\og classique\fg},
g\'en\'eralis\'e par Zhang dans~\cite{zhang95,zhang95b}
\`a des situations singuli\`eres; Abb\`es et Bouche~\cite{abbes-b95}
en ont aussi donn\'e une autre d\'emonstration.
Dans le cas des courbes, Autissier~\cite{autissier2001b}
a consid\'erablement affaibli
les hypoth\`eses d'un tel th\'eor\`eme.

Ainsi, nous allons d'abord traiter le cas d'une m\'etrique ample. 
En utilisant le th\'eor\`eme d'Autissier,
nous allons ensuite \'etudier le cas des courbes:
nous y montrons un th\'eor\`eme d'\'equidistribution des petits points
tr\`es g\'en\'eral; il s'applique aussi bien aux hauteurs d\'efinies
par des syst\`emes dynamiques sur la droite projective,
ou qu'\`a la hauteur de N\'eron--Tate sur une courbe elliptique.
Toutefois, faute d'une description tr\`es explicite des mesures limites,
nous donnerons pour finir un r\'esultat d'\'equidistribution
dans le {\og graphe de r\'eduction\fg} de Chinburg--Rumely.
La situation est d\'ej\`a  int\'eressante dans le cas de courbes
elliptiques \`a r\'eduction semi-stable, dont le graphe de r\'eduction
est un cercle : les points de petite hauteur s'y \'equidistribuent 
(pour la mesure invariante par rotation). Cela signifie
que les r\'eductions 
des points de petite hauteur 
modulo une place de mauvaise r\'eduction 
{\og visitent\fg} r\'eguli\`erement toutes les composantes irr\'eductibles
de la limite inductive des mod\`eles de N\'eron sur 
les corps de nombres contenus dans~$\bar\Q$.

\bigskip

Pour finir cette introduction,
citons deux approches parall\`eles \`a la distribution
des points de petite hauteur pour les syst\`emes dynamiques
sur la droite projective : il s'agit des travaux
de M.~Baker et R.~Rumely~\cite{baker-r2004}
d'une part, 
et de C.~Favre et J.~Rivera-Letelier~\cite{favre-r-t2004}
d'autre part.
Ces deux groupes d'auteurs utilisent une th\'eorie
du potentiel sur la droite projective de Berkovich
qu'ils ont d\'efinie, ce qui leur permet de g\'en\'eraliser
les d\'emonstrations de Bilu
et Baker-Hsia 
(\cite{bilu97}, \cite{baker-h2003}).
Citons enfin un travail en cours d'A.~Thuillier dans
laquelle est d\'evelopp\'ee une th\'eorie du potentiel
sur les courbes alg\'ebriques de genre arbitraire,
dans l'esprit du livre~\cite{rumely1989} de Rumely.
On y trouvera des compl\'ements aux r\'esultats de cet article.

\bigskip

Apr\`es une premi\`ere version, incompl\`ete et partiellement erron\'ee,
c'est au cours d'un s\'ejour \`a Bombay, au sein d'un projet
\textsc{Cefipra} consacr\'e aux groupes alg\'ebriques,
que j'ai de nouveau pu r\'efl\'echir \`a ces questions.
J'ai aussi eu l'occasion d'exposer les r\'esultats de cet
article au cours du premier semestre 2004 dans des conf\'erences \`a 
New-York et M\"unster.
Merci \`a ces institutions pour leur accueil chaleureux.

Merci enfin \`a A.~Thuillier pour une remarque d\'ecisive
qui a permis la d\'emonstration du th.~\ref{theo.courbes},
ainsi qu'\`a P.~Autissier pour ses commentaires et suggestions.

\section{Construction de mesures}

Soit $K$ un corps muni d'une valeur absolue ultram\'etrique
pour laquelle il est complet. On note $K^\circ$
l'anneau des~$a\in K$ tels que $\abs a\leq 1$,
$K^{\circ\circ}$ l'ensemble des~$a\in K$ tels que
$\abs a<1$. C'est un id\'eal maximal de~$K^\circ$ ; on notera
$\tilde K=K^\circ/K^{\circ\circ}$
le  corps r\'esiduel.

Soit $X$ un $K$-espace analytique strict, au sens de
Berkovich~\cite{berkovich1990}.
Les espaces que nous aurons \`a consid\'erer seront toujours
suppos\'es correspondre \`a un espace rigide quasi-compact
et quasi-s\'epar\'e. Puisque le spectre de Berkovich $\mathscr M(\mathscr A)$
d'une $K$-alg\`ebre affino\"{\i}de~$\mathscr A$ est compact,
de tels espaces sont compacts.
On peut de plus les \'etudier via la g\'eom\'etrie formelle
car ils sont la fibre g\'en\'erique
d'un sch\'ema formel admissible quasi-compact sur~$K^\circ$
(voir~\cite{bosch-l1993}, th.~4.1).
Ils seront le plus souvent les analytifi\'es de vari\'et\'es projectives. 
Dans ce cas, tout $K^\circ$-mod\`ele formel est le but d'un
\'eclatement admissible de source un $K^\circ$-sch\'ema projectif
(\cf par exemple~\cite{gubler2003}, prop.~10.5).

Nous supposerons toujours que $K$ admet un sous-corps d\'enombrable
et dense.  C'est le cas pour nos applications \`a la th\'eorie d'Arakelov
puisque $K$ est alors ou bien une extension finie de~$\Q_p$,
ou bien le corps $\C_p$. Dans ces deux cas, la cl\^oture alg\'ebrique
de~$\Q$ dans~$K$ est dense.

Cette hypoth\`ese assure que l'espace~$X$ est m\'etrisable.
Soit en effet $\Omega$ un sous-corps d\'enombrable et dense dans~$K$.
Toute $K$-alg\`ebre strictement affino\"{\i}de $\mathscr A$ est s\'eparable,
c'est-\`a-dire admet une famille dense et d\'enombrable.
Par suite, une semi-norme born\'ee 
   multiplicative sur $\mathscr A$ est d\'etermin\'ee par ses restrictions
   \`a une partie $\{f_n\}$ dense dans la boule unit\'e de~$\mathscr A$,
de sorte
que l'application
   de $\mathscr M(\mathscr A)$ dans~$[0,1]^{\N}$, donn\'ee par
   $x\mapsto \abs{f_n(x)}$ est une injection continue. Comme $\mathscr
M(\mathscr A)$ est compact, cette injection induit un hom\'eomorphisme 
de $\mathscr M(\mathscr A)$ sur son image. Puisque $[0,1]^{\N}$
est m\'etrisable, il s'ensuit que $\mathscr M(\mathscr A)$ est un
espace compact m\'etrisable.
   Le cas d'un $K$-espace strictement analytique compact
s'en d\'eduit car un tel espace poss\`ede un recouvrement 
ferm\'e localement fini par des parties strictement affino\"{\i}des.
Voir aussi~\cite{mainetti2001}.

\bigskip

Les premiers paragraphes sont consacr\'es \`a des rappels
de r\'esultats de Gubler~\cite{gubler1998} et Zhang~\cite{zhang95b}.

\subsection{M\'etriques}
Soit $L$ un faisceau inversible sur~$X$.
Une m\'etrique continue sur~$L$ est la donn\'ee, pour tout
ouvert $U$ de~$X$ et toute section $s\in\Gamma(U,L)$
d'une fonction continue $\norm s\colon U\ra\R_+$ telle
que $\norm{fs}=\abs f\norm s$ si $f\in\mathscr C^0(U,K)$
et qui ne s'annule pas si~$s$ ne s'annule pas.
Le produit tensoriel de deux faisceaux inversibles m\'etris\'es,
le dual d'un faisceau inversible m\'etris\'e sont munis
d'une m\'etrique continue canonique.

Soit $e$ un entier strictement
positif et soit $(\mathfrak X,\mathfrak L)$ un mod\`ele formel de~$(X,L^e)$
sur~$K^\circ$ :
$\mathfrak X$ est un sch\'ema formel admissible sur~$K^\circ$
de fibre g\'en\'erique~$X$
et $\mathfrak L$ est un faisceau inversible sur~$\mathfrak X$
de fibre g\'en\'erique~$L^e$.
Cette donn\'ee munit~$L$ d'une m\'etrique continue pour
laquelle les sections locales de norme~$\leq 1$ sont celles
qui sont {\og enti\`eres\fg} (\cite{gubler1998}, Lemma~7.4) :
si $\mathfrak U$ est un ouvert formel
de~$\mathfrak X$ de fibre g\'en\'erique~$U$, $\sigma\in\Gamma(\mathfrak
U,\mathfrak L)$ une trivialisation de~$\mathfrak L$ sur~$\mathfrak U$,
alors $\norm{\sigma}(x)=1$ pour tout $x\in U$.
Une telle m\'etrique sera dite \emph{formelle}.

Tout faisceau inversible poss\`ede une m\'etrique formelle
(\cite{gubler1998}, Lemma~7.6).

Toute fonction continue $f$ sur~$X$ d\'efinit une m\'etrique
continue sur le faisceau inversible trivial $\mathscr O_X$,
d\'efinie par $\norm{1}(x)=e^{-f(x)}$, et r\'eciproquement.
On notera $a\colon\mathscr C^0(X)\ra\bPic(X)$ l'homomorphisme
de groupes ainsi d\'efini.
Il r\'esulte du th\'eor\`eme 7.12 de~\cite{gubler1998}
que les m\'etriques formelles
sont denses dans l'ensemble des m\'etriques continues
sur un fibr\'e en droites.
(c'est une application du th\'eor\`eme de Stone-Weierstrass sur l'espace
compact~$X$).

\subsection{M\'etriques alg\'ebriques; m\'etriques int\'egrables}
Soit $X$ une $K$-vari\'et\'e projective et notons $X^{\an}$ l'espace
de Berkovich qui lui est associ\'e.
Soit $L$ un fibr\'e en droites sur~$X$.

Parmi les m\'etriques sur~$L$, on appellera \emph{m\'etriques alg\'ebriques}
celles qui sont d\'efinies par un mod\`ele
$(\mathscr X,\mathscr L)$, o\`u $\mathscr X$ est
(le sch\'ema formel associ\'e \`a) un sch\'ema projectif et
$\mathscr L$ un fibr\'e en droites sur $\mathscr X$
tel que $\mathscr L \simeq L^e$, pour un entier $e\geq 1$.

Tout fibr\'e en droite admet une m\'etrique alg\'ebrique;
plus pr\'ecis\'ement, toute m\'etrique formelle est alg\'ebrique.

Suivant Zhang (\cite{zhang95b},  (1.3)),
nous dirons qu'une m\'etrique alg\'ebrique est \emph{semi-positive}
si la r\'eduction~$\mathsf L$ de $\mathfrak L$ 
a un degr\'e positif ou nul sur toute courbe 
de la vari\'et\'e propre~$\mathsf X$, r\'eduction de~$\mathfrak X$.
(Voir aussi \cite{gubler1998}, D\'ef.~7.13.)
Nous dirons plus g\'en\'eralement qu'une m\'etrique continue
est semi-positive si elle est limite uniforme
de m\'etriques alg\'ebriques semi-positives.
Nous dirons enfin qu'une m\'etrique sur~$L$ est \emph{int\'egrable}
si $L$ est isom\'etrique au quotient de deux faisceaux inversibles
m\'etris\'es \`a m\'etrique semi-positive.

Nous noterons $\bPic(X)$ (\resp $\bPicalg(X)$, $\bPicfor(X)$
et $\bPicint(X)$) les groupes de classes d'isomorphismes
de faisceaux inversibles m\'etris\'es (\resp dont la  m\'etrique
est alg\'ebrique, formelle, int\'egrable).
On a ainsi des inclusions
\[ \bPicfor(X) = \bPicalg(X) \subset \bPicint(X) \subset\bPic(X). \]
L'inclusion $\bPicalg(X)\subset\bPicint(X)$ r\'esulte de ce que
sur un $K^\circ$-sch\'ema projectif, tout faisceau inversible
est le quotient de deux faisceaux inversibles amples.

\subsection{Mesures}

Soit encore $X$ une $K$-vari\'et\'e projective de dimension~$d$.
Soit $\bar{L_i}$, $1\leq i\leq d$, des faisceaux
inversibles m\'etris\'es sur~$X$ dont les m\'etriques sont alg\'ebriques.
Elles sont d\'efinies par des fibr\'es en droites $\mathscr L_i$
sur un m\^eme $K^0$-sch\'ema projectif $\mathscr X$, o\`u $\mathscr L_i$
est un mod\`ele de $L_i^{e_i}$, pour un certain entier $e_i\geq 1$.
On peut en outre, et nous le ferons syst\'ematiquement, supposer 
que $\mathscr X$ est normal;
notons $\mathsf X$ sa fibre sp\'eciale. 
L'application de r\'eduction $\pi\colon X^\an\ra\mathsf X$ est surjective
(utiliser la prop.~2.4.4 de~\cite{berkovich1990});
de plus, puisque $\mathscr X$ est suppos\'e normal, l'espace $X^\an$
contient pour chaque composante $\mathsf X_j$ de~$\mathsf X$ un unique
point  $\xi_j$ dont la r\'eduction est le point g\'en\'erique de~$\mathsf X_j$:
la semi-norme $\xi_j$ est d\'efinie par la valuation discr\`ete 
le long de~$\mathsf X_j$
sur les alg\`ebres affino\"{\i}des $\Gamma(\pi^{-1}(\mathsf U), \mathscr O_X)$, 
pour $\mathsf U\subset\mathsf X_j$ affine.
Notons enfin $\mathsf L_i$ la r\'eduction de~$\mathscr L_i$
sur~$\mathsf X$.

Dans le cas alg\'ebrique, nos mesures sont des combinaisons
lin\'eaires de masses de Dirac:
\begin{defi}
On d\'efinit une mesure sur $X^\an$ par la formule
\[  c_1(\bar{L_1})\cdots c_1(\bar{L_d}) =
  \frac{1}{e_1\dots e_d} 
\sum_j  \big( c_1(\mathsf L_1)\cdots c_1(\mathsf L_d)|\mathsf X_j \big)
\nu_j \delta_{\xi_j}  , \]
o\`u $\nu_j$ est la multiplicit\'e de~$\mathsf X_j$
et $\delta_{\xi_j}$ la mesure de Dirac normalis\'ee support\'ee en~$\xi_j$.
\end{defi}
Cette mesure d\'epend de mani\`ere $d$-lin\'eaire sym\'etrique des~$\bar{L_i}$.
Il r\'esulte aussi de la d\'efinition que $c_1(\bar{L_1})\cdots
c_1(\bar{L_d})$ est une mesure positive
si les~$\bar{L_i}$ sont des faisceaux inversibles
m\'etris\'es dont les m\'etriques sont alg\'ebriques semi-positives.
Enfin, la masse totale de cette mesure est \'egale
au degr\'e $(c_1(L_1)\dots c_1(L_d)|X)$.

Plus g\'en\'eralement, si $Z$ est une sous-vari\'et\'e irr\'eductible de~$X$,
de dimension~$e$, on d\'efinit
\[ c_1(\bar L_1)\dots c_1(\bar L_k) \delta_Z 
= i_* c_1(\bar L_1|_Z)\dots c_1(\bar L_e|Z), \]
o\`u $i\colon Z^\an\ra X^\an$ est l'immersion canonique
et $i_*$ est l'application qu'on en d\'eduit sur les mesures.
C'est une mesure de masse totale $(c_1(L_1)\dots c_1(L_k)|Z)$,
positive si les m\'etriques sur les~$\bar L_k$ 
sont alg\'ebriques semi-positives.

\subsection{Lien avec la th\'eorie de l'intersection}

Comme pr\'ec\'edemment,
soit $\mathscr X$ un $K^\circ$-mod\`ele normal de~$X$.
D\'esignons par $Z(\mathscr X)$ le groupe des cycles sur~$\mathscr X$,
\`a coefficients entiers pour les cycles horizontaux
et \`a coefficients r\'eels pour les cycles verticaux.

Soit $\bar L$ un fibr\'e en droite sur~$X$ muni d'une m\'etrique
d\'efinie par un mod\`ele $\mathscr L$ de~$L^e$ sur~$\mathscr X$.
Toute section m\'eromorphe $s$ de~$L$ sur~$X$ a un diviseur
$\hdiv (s)$ qui est un diviseur de Cartier sur~$\mathscr X$.
Une composante $\mathsf X_j$ de la fibre sp\'eciale de~$\mathscr X$,
appara\^{\i}t avec multiplicit\'e $-\log\norm{s}(\xi_j)$,
o\`u $\xi_j$ est l'unique point de~$X$ dont la r\'eduction
est le point g\'en\'erique de~$\mathsf X_j$.

Dans le cas o\`u $K^\circ$ est un anneau de valuation
discr\`ete d'uniformisante~$\pi$, 
la fonction~$\pi$ d\'efinit $\nu_j \mathsf X_j$ au voisinage
du point g\'en\'erique~$\mathsf X_j$; si $s$ est une \'equation
locale de~$\mathsf X_j$ au voisinage de ce point g\'en\'erique,
on a donc $-\log\norm{s}(\xi_j)=-\log\abs\pi /\nu_j$.
Cela montre que la d\'efinition de $\hdiv(s)$ co\"{\i}ncide
avec celle d\'efinie en th\'eorie de l'intersection arithm\'etique
\`a un facteur $\log\abs\pi^{-1}$ pr\`es.
(Apr\`es traduction dans le langage des espaces de Berkovich,
\cf aussi les prop.~7.2 et d\'ef.~7.8 de~\cite{gubler2003}.)

Soit $\bar {L_0},\dots,\bar{L_d}$ des fibr\'es
en droites m\'etris\'es sur~$X$, dont les m\'etriques sont d\'efinies
par des mod\`eles sur~$\mathscr X$.
Pour tout~$i$, soit $s_i$ une section m\'eromorphe non nulle de~$L_i$.
Soit $Z$ une sous-vari\'et\'e de $X$.

Si les $s_i$ s'intersectent proprement sur~$Z$, 
la th\'eorie de l'intersection arithm\'etique permet de d\'efinir
la hauteur locale de~$Z$ relativement \`a ces sections,
par exemple en extirpant de la d\'efinition globale~(2.3.1)
de la hauteur dans~\cite{bost-g-s94} la composante en la place
qui nous int\'eresse. La prop.~2.3.1, formule~(2.3.8) de \emph{loc. cit.}
montre que cette hauteur locale est d\'efinie par r\'ecurrence, par la formule
\[ (\hdiv(s_0)\dots \hdiv(s_d)|Z)
 = (\hdiv(s_1)\dots\hdiv(s_d)|\div(s_0|_Z)) 
       - \int_{X} \log\norm{s_0} c_1(\bar {L_1}) \dots c_1(\bar{L_d})\delta_Z
.\]
Cette derni\`ere formule de r\'ecurrence vaut m\^eme si $K^\circ$
n'est pas un anneau de valuation discr\`ete; sous une forme
essentiellement \'equivalente, elle est \`a la base de la d\'efinition
par W.~Gubler de hauteurs locales de cycles dans~\cite{gubler1998}.

L'int\'er\^et de l'introduction des mesures $c_1(\bar {L_1}) \dots
c_1(\bar{L_d})\delta_Z$ vient de ce que la formule pr\'ec\'edente
vaut aussi dans le cas complexe, r\'einterpr\'et\'ee
de sorte que $\bar L$ d\'esigne un faisceau inversible muni
d'une m\'etrique $\mathscr C^\infty$, $c_1(\bar{L_i})$ est son
courant de courbure et $\delta_Z$ est le courant d'int\'egration sur~$Z$.

Un cas particulier de cette construction justifie une notation
suppl\'ementaire: lorsque les fibr\'es m\'etris\'es $\bar{L_i}$
sont tous \'egaux \`a un m\^eme fibr\'e m\'etris\'e $\bar L$
dont le fibr\'e en droites sous-jacent est big, on posera
\[ \mu_{\bar L}= \dfrac1{c_1(L)^d} c_1(\bar L)^d. \]
C'est une mesure de masse totale~$1$ sur~$X$;
si la m\'etrique de~$\bar L$ est semi-positive, c'est m\^eme
une mesure de probabilit\'e.

\subsection{Passage \`a la limite}

Supposons maintenant que les $\bar{L_i}$ soient
des faisceaux inversibles m\'etris\'es arbitraires sur~$X$. Peut-on
\'etendre la d\'efinition des hauteurs locales et des
mesures donn\'ee ci-dessus? En reprenant les
arguments de Zhang dans~\cite{zhang95b}, (1.4),
nous allons voir que la r\'eponse est essentiellement \emph{oui}
si les m\'etriques sont int\'egrables.

Soit $Z$ une sous-vari\'et\'e irr\'eductible de~$X$ de dimension~$d$.
Soit $s_0,\dots,s_d$ des sections m\'eromorphes de $L_0,\dots,L_d$
dont les diviseurs s'intersectent proprement sur~$Z$.

\begin{prop}
Supposons que pour tout~$i$, $\bar{L_i}$ est muni d'une m\'etrique
semi-positive.
Pour toute suite de m\'etriques alg\'ebriques semi-positives 
$\bar{L_{i,n}}$ sur les~$L_i$ qui converge vers les m\'etriques donn\'ees,  
\begin{enumerate} \def\theenumi{\alph{enumi}}\def\labelenumi{\theenumi)}
\item La suite 
$ (\hdiv(s_{0,n})\cdots \hdiv(s_{d,n})|Z)$
des hauteurs locales de~$Z$
converge vers un nombre r\'eel que l'on notera
\[ (\hdiv(s_0)\cdots \hdiv(s_{d})|Z). \]
\item La suite 
$c_1(\bar {L_{1,n}})\dots c_1(\bar{L_{d,n}})\delta_Z$ 
des mesures sur~$Z$  converge vers une mesure sur~$Z$
que l'on notera $c_1(\bar{L_1})\dots c_1(\bar{L_d})\delta_Z$.
\end{enumerate}
\end{prop}
\begin{proof}
\emph a)
L'espace des m\'etriques continues sur un fibr\'e en droite~$L$
est muni de la distance 
\[ d(\norm{\cdot},\norm{\cdot}') = \log \sup
\big( \frac{\norm{\cdot}'}{\norm{\cdot}},
\frac{\norm{\cdot}}{\norm{\cdot}'} \big) \]
qui en fait un espace m\'etrique complet.
Montrons alors que la hauteur locale est une fonction
uniform\'ement continue sur le sous-espace des m\'etriques semi-positives
alg\'ebriques ; il en r\'esultera que cette hauteur locale s'\'etend
par continuit\'e de mani\`ere unique \`a l'espace des m\'etriques semi-positives.

Notons donc $(\bar L_0,\dots,\bar L_d)$ et $(\bar L'_0,\dots,\bar L'_d)$
deux familles de m\'etriques alg\'ebriques semi-positives. Notons
$d_j$ la distance des m\'etriques $\bar {L_j}$ et $\bar{L'_j}$.
On notera $\hdiv(s_j)$ le diviseur de $s_j$ calcul\'e pour la m\'etrique
de $\bar{L_j}$, et $\hdiv(s'_j)$ celui calcul\'e pour la m\'etrique
de $\bar{L'_j}$.
On a  alors
\begin{multline*}
 (\hdiv(s'_0)\cdots \hdiv(s'_d)|Z)
 - (\hdiv(s_0)\cdots \hdiv(s_d)|Z)
\\
= \sum_{k=0}^d (\hdiv(s'_0)\cdots \hdiv(s'_{k}) \hdiv(s'_{k+1})\dots
\hdiv(s_d)|Z)
 - (\hdiv(s'_0)\cdots \hdiv(s'_k) \hdiv(s_{k+1})\dots
\hdiv(s_d)|Z).
\end{multline*}
D'apr\`es la proposition~9.5 de~\cite{gubler1998}, le terme d'indice~$k$
est major\'e par 
\[ d_j (c_1(L_0)\dots c_1(L_{k-1})c_1(L_{k+1})\dots c_1(L_d)|Z), \]
qui est un multiple constant de la distance~$d_j$.
L'uniforme continuit\'e en r\'esulte, d'o\`u l'existence
des hauteurs locales pour des m\'etriques semi-positives.

\emph b)
Pour montrer qu'une suite de mesure sur~$X^\an$ converge, il suffit
de prouver qu'elle n'a qu'une valeur d'adh\'erence, car $X^\an$
\'etant un espace topologique compact et m\'etrisable, l'espace
des mesures de probabilit\'e sur~$X^\an$ est compact pour la topologie
de la convergence vague.
Soit $\nu$ une valeur d'adh\'erence de la suite $c_1(\bar {L_{1,n}})\dots
c_1(\bar {L_{d,n}})\delta_Z$.

Soit $\bar{L_0}$ et $\bar{L'_0}$ deux m\'etriques alg\'ebriques semi-positives
sur un fibr\'e  en droites tr\`es ample~$L_0$. Soit $s_0$
une section globale de~$L_0$ telle que $s_0,\dots,s_d$
se coupent proprement sur~$Z$.
On peut donc \'ecrire $\hdiv(s_0)'=\hdiv(s_0)+a(f)$, o\`u $f$
est une fonction continue sur~$X^\an$.
Par passage \`a la limite dans les hauteurs locales, on a donc
\begin{align*}
\int_X f\nu &=
 \lim_n \int a(f) c_1(\bar{L_{1,n}})\dots c_1(\bar{L_{d,n}})\delta_Z \\
& = \lim_n 
(\hdiv(s_0)'\hdiv(s_{1,n})\dots\hdiv(s_{d,n})|Z)
-(\hdiv(s_0)\hdiv(s_{1,n})\dots\hdiv(s_{d,n})|Z)  \\
& = (\hdiv(s_0)'\hdiv(s_1)\dots\hdiv(s_d)|Z)
-(\hdiv(s_0)\hdiv(s_1)\dots\hdiv(s_d)|Z) ,
\end{align*}
ce qui montre que $\int f\nu$ est d\'etermin\'e par l'accouplement
de hauteurs locales. Comme l'espace vectoriel engendr\'e par les fonctions~$f$
de cette forme est dense dans l'espace des fonctions continues sur $X^\an$,
il existe au plus une telle mesure~$\nu$, ainsi qu'il fallait d\'emontrer.
\end{proof}

Supposons que pour tout~$j$, $\bar{L_j}$ soit un fibr\'e en droites 
munis d'une m\'etrique int\'egrable
et $s_j$ est une section m\'eromorphe des~$L_j$.
Si $\div(s_0),\dots,\div(s_d)$ s'intersectent proprement
sur~$Z$, on en d\'eduit par multilin\'earit\'e une hauteur locale
$(\hdiv(s_0)\dots\hdiv(s_d)|Z)$,
ainsi qu'une mesure
$c_1(\bar L_1)\dots c_1(\bar L_d)\delta_Z$.

\subsection{Propri\'et\'es}

Soit $f\colon X\ra Y$ un morphisme g\'en\'eriquement fini de $K$-vari\'et\'es
projectives int\`egres.
Soit $\bar{L_i}$ ($0\leq i\leq d$) des faisceaux
inversibles m\'etris\'es sur~$Y$ dont les m\'etriques sont int\'egrables.
On a alors l'\'egalit\'e de mesures sur $Y^\an$:
\[ f_*\big( c_1(f^*\bar{L_1})\cdots c_1(f^*\bar{L_d}) \big) 
 =  \deg(f) c_1(\bar{L_1})\cdots c_1(\bar{L_d}). \]

Supposons que $f$ soit fini et que $X$ et $Y$ soient normale;
on peut alors d\'efinir un homomorphisme trace 
$\mathscr C^0(X^\an)\ra\mathscr C^0(Y^\an)$.
Soit $f'\colon X'\ra Y$ la normalisation de~$Y$ dans une cl\^oture
galoisienne de l'extension finie $K(Y)\subset K(X)$. 
Le morphisme~$f'$ se factorise en $f\circ g$, o\`u $g\colon X'\ra X$
est un morphisme fini. En outre, le groupe $\Gal(K(X')/K(Y))$
agit sur~$X'$ par $Y$-morphismes; il agit donc aussi
sur $(X')^\an$. Si $\phi\in\mathscr C^0(X^\an)$, la fonction
\[  \tilde\phi= \frac1{[K(X'):K(X)]}\sum_{\gamma\in\Gal(K(X')/K(Y))}
           \phi\circ\gamma \]
sur $(X')^\an$ est continue et descend \`a~$Y^\an$; le morphisme
$(X')^\an\ra Y^\an$ est fini, donc propre (prop.~3.4.7
et corollaire~3.3.8 de~\cite{berkovich1990}) si bien que l'unique fonction
$f_*(\phi)$ sur~$Y^\an$ telle que $\tilde\phi=f_*(\phi)\circ f'$
est continue. 
Par dualit\'e, on en d\'eduit un homomorphisme $f^*$ sur les mesures.

On a alors
\[  f^*(c_1(\bar{L_1})\cdots c_1(\bar{L_d}))
= c_1(f^*\bar{L_1})\cdots c_1(f^*\bar{L_d}) .
 \]
Ces formules sont en effet v\'erifi\'ees si les m\'etriques sont alg\'ebriques,
ainsi qu'il r\'esulte de la formule de projection.
Le cas g\'en\'eral s'en d\'eduit par passage \`a la limite.
              
La compatibilit\'e de ces mesures au produit est un peu plus subtile
(mais gu\`ere plus difficile).
Notons $\abs X$ l'espace topologique sous-jacent \`a un espace
de Berkovich~$X$. Si $X$ et $Y$ sont deux espaces de Berkovich,
on dispose alors d'une application canonique, 
surjective et propre, $\alpha\colon \abs{X\times Y}\ra \abs X\times\abs Y$.
Notons $d=\dim X$, $e=\dim Y$, et supposons que ce soient les
analytifi\'es de $K$-vari\'et\'es projectives.
Soit $\bar L$ et $\bar M$ des faisceaux inversibles 
sur $X$ et $Y$ respectivement, munis de m\'etriques
int\'egrables. Notons $\bar L\boxtimes\bar M$
le faisceau inversible qui s'en d\'eduit sur $X\times Y$;
il est muni d'une m\'etrique int\'egrable naturelle et l'on a,
par passage \`a la limite du cas alg\'ebrique, la formule
\[ \alpha_* \big( c_1(\bar L\boxtimes\bar M)^{d+e} \big)
 = \binom{d+e}d 
      c_1(\bar L)^d \otimes c_1(\bar M)^d, \]
o\`u le produit tensoriel du membre de droite
repr\'esente la mesure sur $\abs X\times \abs Y$,
produit tensoriel des mesures $c_1(\bar L)^d$ sur~$X$
et $c_1(\bar M)^e$ sur~$Y$.

\subsection{Comparaison avec l'intersection d'Arakelov}

Il s'agit en fait de la g\'en\'eralisation par Zhang
dans~\cite{zhang95b}
de l'accouplement d'Arakelov qu'ont d\'efini Gillet et Soul\'e
(\cite{gillet-s90,bost-g-s94}).
Dans ce contexte, $F$ est un corps de nombres,
$X$ est une $F$-vari\'et\'e projective,
les~$\bar{L_i}$ des faisceaux inversibles munis de m\'etriques
ad\'eliques int\'egrables (\cf~\cite{zhang95b}).

Une telle m\'etrique sur un fibr\'e inversible $L$ sur~$X$ est la donn\'ee, 
pour toute place $v$ de~$F$,
d'une m\'etrique int\'egrable sur le fibr\'e en droites
sur la vari\'et\'e analytique $X_v$
d\'eduit de~$L$. La condition pour la m\'etrique
d'\^etre ad\'elique signifie
qu'il existe un mod\`ele de $X$ et $L$
sur un ouvert de $\Spec\mathfrak o_F$ d\'efinissant simultan\'ement
presque toutes ces m\'etriques.
En particulier, la m\'etrique $v$-adique de~$L$ est, pour presque
toute place de~$F$, alg\'ebrique.

Pour tout~$i$, soit $s_i$ une section m\'eromorphe de~$L_i$, 
choisies de sorte que les $\div(s_i)$
s'intersectent proprement sur la sous-vari\'et\'e~$Z$ de~$X$.
Pour toute place $v$, la th\'eorie des hauteurs locales (\cf\cite{gubler1997}) 
d\'efinit alors le nombre r\'eel
$\big(\hdiv(s_0)\cdots \hdiv(s_d) |Z\big)_v$;
ces nombres sont presque tous nuls et l'on a la formule
\begin{equation}
\label{eq.localglobal}
 (\hc_1(\bar{L_0})\cdots \hc_1(\bar L_d)|Z)
       = \sum_v \big(\hdiv(s_0)\cdots \hdiv(s_d) |Z\big)_v \log N_v , 
\end{equation}
o\`u $N_v$ est le cardinal du corps r\'esiduel de~$F$ en~$v$.

Dans une suite d'approximations alg\'ebriques amples d'une m\'etrique 
semi-positive donn\'ee, le membre de gauche converge,
ainsi que l'a montr\'e Zhang; la limite est not\'ee
naturellement $ (\hc_1(\bar{L_0})\cdots \hc_1(\bar L_d)|Z)$.
Par ailleurs, les expressions d\'efinissant les hauteurs
locales convergent place par place vers 
la quantit\'e $\big(\hdiv(s_0)\cdots \hdiv(s_d) |Z\big)_v$
d\'efinie plus haut  si $v$ est non archim\'edienne; une variante
de ce raisonnement \'etablirait le cas o\`u $v$ est archim\'edienne.
La formule~\eqref{eq.localglobal}
est donc v\'erifi\'ee pour des m\'etriques ad\'eliques semi-positives.
Sa validit\'e dans le cas de m\'etriques int\'egrables en r\'esulte
par multilin\'earit\'e.

\section{\'Equidistribution : le cas d'une m\'etrique ample}

\begin{theo}
\label{theo.berk-ample}
Soit $X$ une vari\'et\'e projective sur un corps de nombres~$F$.
Soit $L$ un fibr\'e en droites ample sur~$X$ muni d'une m\'etrique
ad\'elique semi-positive.

Soit $(x_n)$ une suite g\'en\'erique de points de~$X(\bar F)$
telle que $h_{\bar L}(x_n)\ra h_{\bar L}(X)$.
Alors, pour toute place~$v$ de~$F$ en laquelle la m\'etrique de~$L$ 
est \emph{ample}, la suite des mesures $(\mu_{x_n})$ sur
l'espace $X_v$ converge faiblement vers
la mesure de probabilit\'e $\mu_{\bar L_v}$. 
\end{theo}

Dans le cas d'une place archim\'edienne, ce th\'eor\`eme
est d\^u \`a Szpiro, Ullmo et Zhang (\cf\cite{szpiro-u-z97}).
Nous nous contenterons donc de donner la d\'emonstration
dans le cas non archim\'edien.

Rappelons tout d'abord
ce que vaut la mesure limite, en supposant que $v$ est une place ultram\'etrique.
Soit $(\mathscr X,\mathscr L)$ un mod\`ele sur $\mathfrak o_{F_v}$ d\'efinissant
la m\'etrique $v$-adique de $\bar L$.
Notons $\mathsf Y_1,\dots,\mathsf Y_r$ les composantes irr\'eductibles
de sa fibre sp\'eciale en~$v$, $m_1,\dots,m_r$ leurs multiplicit\'es.
Pour tout~$i\in\{1,\dots,r\}$, $\mathsf Y_i$ est la r\'eduction
d'un unique point, $\eta_i$, de l'espace de Berkovich $X_v$.
Alors,
\[ \mu_{\bar L} = \frac{c_1(\overline L)^{\dim X}}{c_1(L)^{\dim X}} = \sum_{i=1}^r m_i \frac{(c_1(\mathscr L)^{\dim \mathsf Y_i}|\mathsf Y_i)}{c_1(L)^{\dim X}|X)} \delta_{\eta_i}.
\]

\begin{exem}
Supposons que $X=\P^d$ et que $\bar L=\mathscr O(1)$,
muni de sa m\'etrique $L^\infty$ qui donne lieu \`a la hauteur de Weil
usuelle. Cette m\'etrique est semi-positive et ample aux places
ultram\'etriques, car d\'efinie par le mod\`ele $(\P^d_{\mathfrak o_F},\mathscr
O(1))$. Pour toute place finie~$v$,
la mesure limite $\mu_{\bar L}$ sur l'espace
de Berkovich~$\P^d_{F_v}$ est la masse de Dirac au point canonique de 
$\P^d_{F_v}$ dont la r\'eduction en~$v$ est le point g\'en\'erique
de~$\P^d_{k_v}$. Nous noterons ce point~$G_v$ car lorsque~$d=1$,
il correspond \`a la norme de Gauss sur l'anneau des s\'eries formelles
restreintes $F_v\{T\}$.

On sait que $h_{\bar L}(\P^d)=0$ (\cf\cite{zhang92} ou \cite{maillot2000}
pour le cas plus g\'en\'eral des vari\'et\'es toriques auquel, du reste,
cet argument s'applique aussi).
Il d\'ecoule du th\'eor\`eme que pour toute suite g\'en\'erique $(x_n)$ de points de
$\P^d(F_v)$ dont la hauteur tend vers~$0$
et toute place finie~$v$, la suite
de mesure $(\mu_{x_n})$ sur l'espace de Berkovich $\P^d_{F_v}$
converge vers la mesure de Dirac au point~$G_v$.
C'est l'analogue ultram\'etrique d'un th\'eor\`eme de Bilu.

Joint \`a la solution (Laurent, Bilu) 
de la {\og conjecture de Bogomolov sur les tores\fg}, notre
th\'eor\`eme s'\'etend aux suites ({\og strictes\fg}) de points de~$\gm^d$
dont la hauteur tend vers~$0$ dont aucune sous-suite n'est
contenue dans un sous-groupe alg\'ebrique strict. 
En effet, cette conjecture affirme que de telles suites
sont g\'en\'eriques.
Notons enfin que si le th\'eor\`eme de Bilu en fournit une d\'emonstration,
il ne semble pas que ce soit le cas pour notre th\'eor\`eme.
\end{exem}

\begin{exem}
On peut aussi appliquer notre th\'eor\`eme au cas d'une vari\'et\'e
ab\'elienne, \emph{en toute place finie~$v$ o\`u elle a bonne r\'eduction}
(quitte \`a effectuer d'abord une extension
des scalaires, bonne r\'eduction potentielle suffit).
La mesure limite est alors concentr\'ee au point de l'espace
de Berkovich dont la r\'eduction est le point g\'en\'erique du mod\`ele
de la fibre sp\'eciale.
Les m\^emes remarques concernant la conjecture de Bogomolov
s'appliquent, ladite conjecture \'etant dans ce cas un th\'eor\`eme
de S.~Zhang (\cite{zhang98}, voir aussi~\cite{david-p98}).
\end{exem}

La d\'emonstration de ce th\'eor\`eme suit l'approche des th\'eor\`emes
d'\'equidistribution
initi\'ee par Szpiro, Ullmo et Zhang et passe par plusieurs \'etapes.

Comme l'espace de Berkovich~$X_{F_v}$  est un espace
topologique compact et m\'etrisable, l'ensemble des mesures de probabilit\'e
sur $X_{F_v}$  est compact pour la topologie vague
et il suffit 
de prouver que  la mesure $\mu_{\bar L}$
est l'unique valeur d'adh\'erence de 
la suite $(\mu_{x_n})$. Pour cela, on peut en outre supposer
que la suite $(\mu_{x_n})$ converge vaguement vers une mesure $\mu$,
il faut alors d\'emontrer que $\mu=\mu_{\bar L}$.
Le th\'eor\`eme~\ref{theo.berk-ample} r\'esulte alors des trois
lemmes suivants.

\medskip

Soit $\mathsf V$ un ferm\'e de la fibre sp\'eciale $\mathsf X=\mathscr X\otimes k_v$
de~$\mathscr X$ en~$v$;
nous appellerons tube de~$\mathsf V$, et noterons $\tube{\mathsf V}$,
l'image r\'eciproque de~$\mathsf V$ par l'application de r\'eduction
$X_{F_v}\ra \mathsf X$.

\begin{lemm}
Si $\mathsf V$ ne contient pas de composante irr\'eductible
de~$\mathsf X$, $\mu(\tube{\mathsf V})=0$.
\end{lemm}
(On dira qu'un tel tube est \emph{petit} ;
cette notion d\'epend du choix de~$\mathscr X$.)
\begin{proof}
Soit $\mathscr X'$ l'\'eclatement de~$\mathscr X$ le long de~$\mathsf V$
et soit $\mathsf E$ le diviseur exceptionnel. Le fibr\'e
inversible $\mathscr O_{\mathscr X'}(-\mathsf E)$ d\'efinit
une m\'etrique $v$-adique sur le fibr\'e trivial $\mathscr O_X$;
nous noterons $\phi_{\mathsf V}$ la norme de la section~$1$.
On a $\phi_{\mathsf V}\geq 1$, et l'in\'egalit\'e stricte
vaut pr\'ecis\'ement aux points dont la r\'eduction est dans~$\mathsf V$.
Pour tout nombre r\'eel~$a$,
notons $\bar L_a$ le faisceau~$L$ muni de la m\'etrique ad\'elique
de~$L$, sauf en~$v$ o\`u on la multiplie par~$\phi_{\mathsf V}^a$.

La description
de $\mathscr X'$ comme le spectre projectif d'une alg\`ebre gradu\'ee
entra\^{\i}ne que $\mathscr O_{X'}(-\mathsf E)$ est relativement ample
pour le morphisme d'\'eclatement $\pi\colon\mathscr X'\ra\mathscr X$.
Par suite, pour $n$ assez grand, $\pi^*\mathscr L^{\otimes n}(-\mathsf E)$
est ample sur~$\mathscr X'$. En particulier,
pour tout nombre rationnel $a>0$ assez petit,
$\bar L_a$ est un fibr\'e inversible muni d'une m\'etrique
ad\'elique semi-positive.

Le th\'eor\`eme de Hilbert-Samuel arithm\'etique, sous la forme qui
lui donne Zhang, s'\'ecrit alors
\[ \liminf_n \left( h_{\bar L}(x_n) - a \mu_{x_n}(\log\phi_{\mathsf V}) \right)
      \geq h_{\bar L_a}(X) .
\]
Par d\'efinition,
\begin{align*}
 \hc_1(\bar L_a)^{\dim X} &= 
\big((\pi^*(\hc_1(\mathscr L))- a\hc_1(\mathscr O(\mathsf E)))^{1+\dim X}
      \big|\mathscr X'\big)  \\
&= \sum_{j=0}^{1+\dim X} (-a)^j \binom{1+\dim X}j 
\big( \pi^*\hc_1(\mathscr L)^{1+\dim X-k}\hc_1(\mathscr O(\mathsf E))^k
       \big| \mathscr X' \big) \\
&= \big( \hc_1(\mathscr L)^{1+\dim X} \big|\mathscr X \big) 
\end{align*}
d'apr\`es la formule de projection appliqu\'ee \`a~$\pi$,
car $\pi_*\hc_1(\mathscr O(\mathsf E))=0$.
Par suite, $h_{\bar L_a}(X)=h_{\bar L}(X)$.
Par hypoth\`ese, $h_{\bar L}(x_n)$ converge vers $h_{\bar L}(X)$
lorsque $n$ tend vers~$+\infty$.
Faisant tendre $a>0$ vers~$0$, on trouve
\[ \mu(\log\phi_{\mathsf V}) \leq 0. \] 
En particulier, l'ensemble des points de~$X$ o\`u $\log\phi_{\mathsf V}>0$
est $\mu$-n\'egligeable. On a donc $\mu(\tube{\mathsf V})=0$, comme il
fallait
d\'emontrer.
\end{proof}

\begin{lemm}
Pour $1\leq i\leq r$, $\mu(\tube{\mathsf X_i})$
est donn\'e par la formule
\[ \mu(\tube{\mathsf X_i}) =
 m_i (c_1(\mathscr L)^{\dim X}|\mathsf X_i)/(c_1(\mathscr L)^{\dim X}|X).
\]
\end{lemm}
\begin{proof}
Fixons un entier $i\in\{1,\dots,r\}$ ; le faisceau inversible
$\mathscr O_{\mathscr X}(-\mathsf X_i)$ est un mod\`ele du faisceau
trivial $\mathscr O_X$ et le munit d'une m\'etrique $v$-adique.
Notons $\phi_i$ la norme de sa section canonique~$1$.
Comme dans la d\'emonstration du lemme pr\'ec\'edent, $\phi_i\geq 1$
et l'on a $\phi_i(x)>1$ si et seulement si $x\in\tube{\mathsf X_i}$.
Pour $a\in\R$, notons alors $\bar L_a$ le fibr\'e en droites $L$
muni de la m\'etrique ad\'elique de $\bar L$, sauf en~$v$
o\`u elle est multipli\'ee par $\phi_i^a$.
La m\'etrique en~$v$ de~$\bar L$ est d\'efinie par un faisceau ample
$\mathscr L$ sur~$\mathscr X$ ; pour $t\in\Q$, $t=p/q$ avec
$p$ et $q\in\Z$, la m\'etrique de $\bar L_t$
est d\'efinie par le faisceau $\mathscr L^q(-p\mathsf X_i)$.
Par suite, $\bar L_t$ est semi-positif pour $\abs t$ assez petit.

Si $x\in X(\bar F)$, on a $h_{\bar L_t}(x)=h_{\bar L}(x)- t \mu_x(
\log\phi_i) $, tandis que
\begin{align*} \hc_1(\bar L_t)^{1+\dim X}
& =  \big(\hc_1(\bar L) - t \hc_1(\mathscr O(\mathsf
X_i)))^{1+\dim X}\big|\mathscr X \big) \\
&= \big(\hc_1(\bar L)^{1+\dim X}\big| \mathscr X\big)  
- t \big(c_1(\mathscr L)^{\dim X}c_1(\mathscr O(\mathsf X_i)) \big|\mathscr
X\big) 
 + \mathrm O(t^2)  \\
&= \hc_1(\bar L)^{1+\dim X} - t m_i (c_1(\mathscr L)^{\dim X}|\mathsf X_i)
 + \mathrm O(t^2). 
\end{align*}
Par suite,
\[ h_{\bar L_t}(X) = h_{\bar L}(X) - t  \frac{m_i (c_1(\mathscr L)^{\dim
X}|\mathsf X_i)}{(c_1(L)^{\dim X}|X)} + \mathrm O(t^2). \]
L'in\'egalit\'e de Zhang (Hilbert-Samuel arithm\'etique),
appliqu\'ee \`a la suite $(x_n)$ et au fibr\'e $\bar L_t$,
pour $t$ assez petit, entra\^{\i}ne alors
\[ \mu(\log\phi_i) = 
         m_i \frac{ (c_1(\mathscr L)^{\dim X}|\mathsf X_i)}
                  { (c_1(L)^{\dim X}|X) }. \]

\end{proof}
\begin{lemm}
Soit $K$ un corps complet pour une valuation ultram\'etrique
non triviale, soit $X$ un $K$-espace de Berkovich,
fibre g\'en\'erique d'un $R$-sch\'ema formel, plat et de type fini, $\mathscr X$.
Notons $\mathsf X_i$ les composantes irr\'eductibles de la fibre
sp\'eciale de~$\mathscr X$ et, pour tout~$i$, soit $\xi_i$ 
le point de~$X$ dont~$\mathsf X_i$ est la r\'eduction.
Soit $\nu$ une mesure sur~$X$ telle que tout petit tube est
n\'egligeable. Alors 
\[ \nu = \sum_i \nu(\tube{\mathsf X_i}) \delta_{\xi_i}. \]
\end{lemm}
\begin{proof}
Les tubes sur l'intersection de deux composantes de la fibre
sp\'eciale sont petits, donc n\'egligeables. Cela permet de supposer
que $X$ est le spectre de Berkovich d'une $K$-alg\`ebre strictement
affino\"{\i}de~$\mathscr A$ et que $\mathscr X=\Spf \mathscr A^\circ$,
la fibre sp\'eciale~$\mathsf X$ de~$\mathscr X$ \'etant en outre irr\'eductible.
Soit $\xi$ le point de~$X$ correspondant et soit $U$ un voisinage
de~$\xi$ dans~$X$. Par d\'efinition de la topologie de~$X$, il existe
un ensemble fini de fonctions analytiques $f_j\in \mathscr A$,
pour $j\in J$,
et un nombre r\'eel $\eps>0$ tels que $U$ contienne l'intersection
\[ \bigcap_{j\in J} \{ x\in X\sozat \abs{f_j}(\xi)-\eps <
\abs{f_j}(x)<\abs{f_j}(\xi)+\eps\}. \]
En fait, $\abs{f_j}(x)\leq\abs{f_j}(\xi)$ pour tout $x\in X$
(\cite{berkovich1990}, prop.~2.4.4; $\{\xi\}$ est la fronti\`ere
de Shilov de~$X$). Par suite, le compl\'ementaire de~$U$ est contenu dans
la r\'eunion finie
\[ \bigcup_{j\in J} \{x\in X \sozat \abs{f_j}(x)<\abs{f_j}(\xi) \}. \]
Quitte \`a remplacer $f_j$ par $af_j^{n_j}$ pour  un certain $n_j\in\N^*$
et un certain $a\in K^*$, on peut en outre supposer que $\abs{f_j}(\xi)=1$.
Alors, la r\'eduction de $f_j$ est une fonction non nulle sur le sch\'ema
irr\'eductible $\mathsf X$ et la r\'eunion ci-dessus est une r\'eunion finie
de petits tubes, donc est n\'egligeable. Ainsi, $\nu(\complement U)=0$.

Cela montre que le support de~$\nu$ est r\'eduit au point~$\xi$,
donc $\nu$ est un multiple de la mesure de Dirac en~$\xi$;
la constante de proportionnalit\'e est \'egale \`a $\nu(X)$.
\end{proof}

\section{\'Equidistribution : le cas des courbes}

Les r\'esultats de ce paragraphe ont \'et\'e provoqu\'es par une
remarque d'Amaury Thuillier selon laquelle, sur les
courbes, le th\'eor\`eme de Hilbert-Samuel arithm\'etique 
vaut sans autre condition de positivit\'e  que l'amplitude g\'en\'erique.
C'est en effet ce qui est d\'emontr\'e par Faltings dans~\cite{faltings84}
lorsque les m\'etriques archim\'ediennes sont admissibles au sens
d'Arakelov,  et g\'en\'eralis\'e
par Autissier dans~\cite{autissier2001b} au cas des m\'etriques
archim\'ediennes continues \`a courbure mesure.

Par un passage \`a la limite, on d\'eduit du r\'esultat
d'Autissier la proposition suivante.
\begin{prop}[In\'egalit\'e fondamentale] \label{prop.autissier}
Soit $X$ une courbe projective lisse sur un corps de nombres~$F$.
Soit $L$ un fibr\'e en droites ample sur~$X$, muni d'une m\'etrique
ad\'elique int\'egrable.
Pour toute suite $(x_n)$ de points distincts de~$X(\bar F)$,
on a 
\[ \liminf h_{\bar L}(x_n) \geq  h_{\bar L}(X). \]
\end{prop}
\begin{proof}
Par d\'efinition, on peut \'ecrire $\bar L=\bar M\otimes \bar N^{-1}$,
o\`u $\bar M$ et $\bar M$ sont deux fibr\'es en droites amples
sur~$X$, munis de m\'etriques ad\'eliques semi-positives.
Soit $(\bar {M_k})$ et $(\bar {N_k})$ des suites de m\'etriques
ad\'eliques alg\'ebriques semi-positives sur $M$ et $N$
qui convergent uniform\'ement vers $\bar M$ et $\bar N$ respectivement.
Pour tout~$k$, posons $\bar L_k=\bar{M_k}\otimes \bar {N_k}^{-1}$.
D'apr\`es la prop.~3.3.3 de~\cite{autissier2001b}, on a, pour tout
entier~$k$,
\[ \liminf h_{\bar {L_k}} (x_n) \geq h_{\bar {L_k}}(X). \]
Lorsque $k$ tend vers l'infini, $h_{\bar {L_k}}(x)$
converge vers $h_{\bar L}(x)$, uniform\'ement en $x\in X(\bar F)$,
et $h_{\bar{L_k}}(X)$ converge vers $h_{\bar L}(X)$, par lin\'earit\'e
du produit d'intersection de Zhang et convergence dans le cas
semi-positif.
Il en r\'esulte que
\[ \liminf h_{\bar {L}} (x_n) \geq  h_{\bar {L}}(X),
\]
ce qui est l'in\'egalit\'e de la proposition.
\end{proof}

\begin{theo}\label{theo.courbes}
Soit $X$ une courbe projective lisse sur un corps de nombres~$F$.
Soit $L$ un fibr\'e en droites ample sur~$X$, muni d'une m\'etrique
ad\'elique int\'egrable.
Pour toute suite $(x_n)$ de points distincts de $X(\bar F)$
telle que  $h_{\bar L}(x_n)$ converge vers $h_{\bar L}(X)/2\deg L$,
et toute place~$v$ de~$F$,
la suite de mesures $(\mu_{x_n})$ sur l'espace $X_v$
converge vaguement vers la mesure $c_1(\bar L)/\deg L$.
\end{theo}
\begin{proof}
La d\'emonstration reprend la m\'ethode inaugur\'ee par Szpiro, Ullmo et Zhang.
Les fonctions continues sur $X_v$ d\'efinies par des m\'etriques
alg\'ebriques sont denses dans l'ensemble des fonctions continues.
Il suffit donc de prouver que pour toute telle fonction {\og
alg\'ebrique\fg}~$f$,
 $\mu_{x_n}(f)$ converge vers $\int_{X_v} f c_1(\bar L)/\deg L$.
Appliquons l'in\'egalit\'e de la proposition pr\'ec\'edente 
au fibr\'e m\'etris\'e $\bar L \otimes a(tf)$, pour $t\in\R$.
On obtient
\[ \liminf_n \big( h_{\bar L}(x_n)  + t\mu_{x_n}(f) \big)
 \geq \frac1{2\deg L} \big( h_{\bar L}(X) + 2t \int_{X_v} f c_1(\bar L)
       + \text{terme en~$t^2$.} \]
Par suite,
\[ \liminf t \mu_{x_n}(f) \geq t \int_{X_v} f (c_1(\bar L)/\deg L)
 + \mathrm O(t^2) .\]
Il reste \`a faire tendre~$t$ vers~$0$ d'abord
par valeurs inf\'erieures, puis par valeurs sup\'erieures.
\end{proof}

Remarquons que la m\^eme d\'emonstration fournit un r\'esultat ad\'elique,
\`a savoir l'\'equidistribution d'une suite de petits points dans
le produit des espaces $X_v$ en un nombre fini de places
distinctes de~$F$.

\subsection{Syst\`emes dynamiques}

Cela s'applique notamment aux hauteurs normalis\'ees par
un syst\`eme dynamique alg\'ebrique sur~$\P^1$.
On obtient ainsi l'\'equidistribution des points de petite hauteur,
aux places ultram\'etriques, sous
des hypoth\`eses moins restrictives que  les th\'eor\`emes de Baker-Hsia
dans~\cite{baker-h2003}.
Bien s\^ur, l'\'equidistribution archim\'edienne r\'esultait d\'ej\`a
du th\'eor\`eme d'Autissier (\cite{autissier2001b}, prop.~3.3.3).

Pour pouvoir appliquer ce th\'eor\`eme sur des courbes plus g\'en\'erales
que la droite projective, il importe de disposer de hauteurs
normalis\'ees int\'eressantes. On en construit en consid\'erant
des \emph{correspondances}. Soit $X$ un courbe
projective lisse et  soit $(p_1,p_2)\colon T\ra X\times X$
une correspondance de bidegr\'e~$(a,b)$. Si $L$ est un fibr\'e inversible
sur~$X$, le fibr\'e $T^*L$ est d\'efini comme $N_{p_1}(p_2^*L)$,
o\`u $N_{p_1}$ d\'esigne la norme (sur la premi\`ere classe de Chern,
c'est $(p_1)_*$).  Si $\bar L$ est un fibr\'e muni d'une m\'etrique alg\'ebrique
(\resp continue), $T^*L$ est muni d'une m\'etrique alg\'ebrique 
(\resp continue); on notera $T^*\bar L$
le fibr\'e m\'etris\'e correspondant. Si $\bar L$ est muni d'une m\'etrique
semi-positive, on v\'erifie ais\'ement que le fibr\'e $T^*\bar L$ est encore
semi-positif; par suite, l'op\'eration $\bar L\ra T^*\bar L$
pr\'eserve les fibr\'es en droites munis de m\'etriques ad\'eliques int\'egrables.

Soit $L$ un fibr\'e inversible sur~$X$ tel que $T^*L$ soit isomorphe \`a~$L^b$.
Soit $\alpha\colon L^b\ra T^*L$ un tel isomorphisme.
\emph{Supposons $b>a$.} il existe alors une unique m\'etrique ad\'elique int\'egrable 
(semi-positive si $L$ est ample) sur~$L$ telle que l'homomorphisme
$\alpha\colon \bar L^b\ra T^*\bar L$ soit une isom\'etrie.
Cette m\'etrique se d\'efinit de mani\`ere analogue \`a la m\'etrique canonique
d'un syst\`eme dynamique, comme la limite de la suite de m\'etriques
$\norm{\cdot}_n$, o\`u $\norm{\cdot}_{n+1}=(\alpha^*T^*\norm{\cdot}_n)^{1/b}$,
la m\'etrique $\norm{\cdot}_0$ \'etant une m\'etrique alg\'ebrique arbitraire,
semi-positive si $L$ est ample.
Il en r\'esulte une hauteur canonique $h_{\bar L}$
dans ce contexte, telle que $h_{\bar L}(T_*x)=bh_{\bar L}(x)$
et $h_{\bar L}(X)=0$; le th\'eor\`eme~\ref{theo.courbes} entra\^{\i}ne
donc un r\'esultat d'\'equidistribution sur les suites de petits points
pour cette hauteur normalis\'ee.

\section{Le graphe de r\'eduction}

\subsection{}
Soit $X$ une courbe projective lisse sur un corps de nombres~$F$.
La mauvaise r\'eduction de~$X$ en une place finie~$v$ peut \^etre
cod\'ee de mani\`ere commode dans son graphe de r\'eduction, tel
que l'ont d\'efini Chinburg et Rumely.
Supposons (quitte \`a \'etendre les scalaires)
que $X$ ait r\'eduction semi-stable en~$v$ et notons~$p$
la caract\'eristique du corps r\'esiduel en~$v$.
Le graphe de r\'eduction est un graphe m\'etris\'e dont les sommets
correspondent aux composantes irr\'eductibles de la fibre
sp\'eciale du mod\`ele minimal r\'egulier (semi-stable), deux sommets
sont reli\'es par autant d'ar\^etes de longueur~$1/v(p)$
que les composantes correspondantes ont de points d'intersection communs.
La r\'esolution explicite des singularit\'es de la forme $xy=a$,
$a\in\mathfrak o_F$, montre que la construction du graphe de r\'eduction
commute aux extensions finies.
(Pr\'ecis\'ement, apr\`es une extension finie d'indice de ramification~$e$,
on doit \'eclater les points d'intersection $\lfloor e/2\rfloor$ fois,
cela remplace un segment de longueur~$\ell$ dans le graphe par
$e$ segments de longueur~$\ell/e$.)

Dans son article~\cite{zhang93}, S.~Zhang incorpore \`a
la g\'eom\'etrie d'Arakelov de l'analyse sur ce graphe,
ce qui permet de prendre en compte des m\'etriques ad\'eliques
naturelles, telles la hauteur de N\'eron--Tate d'une courbe
elliptique lorsque celle-ci n'a pas bonne r\'eduction partout.
Du point de vue des fibr\'es en droites m\'etris\'es, on
peut g\'en\'eraliser ce point de vue l\'eg\`erement,
en consid\'erant des mod\`eles r\'eguliers dont la fibre
sp\'eciale ait des croisements normaux stricts: un tel mod\`ele
permet de d\'efinir un graphe de r\'eduction sur lequel
la th\'eorie d'Arakelov non archim\'edienne  dans
le style de celle Zhang peut s'appliquer: c'est d'ailleurs
ce que fait Rumely dans~\cite{rumely1995}.

\subsection{}
Fixons une place ultram\'etrique $v$ de~$F$
et soit $R$ le graphe de r\'eduction de~$X$ \`a la place~$v$.

Soit $E$ une extension finie de~$F$ et soit $\mathscr X_E$
le mod\`ele semi-stable minimal de~$X_E$ sur~$\mathfrak o_E$.
L'application de sp\'ecialisation  
$\sp\colon X(E) \ra R$ est d\'efinie de la fa\c{c}on suivante:
un point $x\in X(E)$ se sp\'ecialise sur une composante bien
d\'efinie de la fibre sp\'eciale de~$\mathscr X_E$ en~$v$,
et $\sp(x)$ est le point correspondant du graphe. Cette
d\'efinition s'\'etend par lin\'earit\'e en un homomorphisme
toujours not\'ee $\sp$, du groupe des diviseurs sur~$X_{\bar F}$
vers le groupe des diviseurs sur le graphe; cet homomorphisme
pr\'eserve le degr\'e. Si $D=\sum n_i P_i$ est un diviseur
sur~$X_{\bar F}$, on notera $\mu_D$ la somme des mesures de Dirac 
$\sum n_i\delta_{\sp(P_i)}$ sur~$R$; c'est une mesure positive.
\`A un point alg\'ebrique~$x\in X(\bar F)$, vu comme un diviseur
de degr\'e~$[F(x):F]$ sur~$X_{\bar F}$, nous attachons la mesure
de probabilit\'e $([F(x):F])^{-1}\mu_x$. 

Suivant~\cite{zhang93}, le graphe de r\'eduction permet de d\'ecrire
certaines classes de m\'etriques ad\'eliques int\'eressantes.
Fixons un fibr\'e en droites~$L$ sur~$X$.

Tout fibr\'e en droites~$\mathscr L$ sur~$\mathscr X_E$
qui \'etend une puissance~$L^{\otimes n_E}$ de~$L$
d\'efinit une m\'etrique sur~$L$. Pour cette m\'etrique, les sections
de norme~$\leq 1$ en un point $P\in X(E)$, 
sont les suivantes: \'etendons d'abord $P$ en un morphisme de~$\mathfrak
o_{E}$-sch\'emas $\eps_P\colon \Spec\mathfrak o_{E'}\ra\mathcal X_E$,
une section~$s$ est de norme~$\leq 1$ en~$P$ si et seulement si
$\eps_P^* s^{n_E}$ s'\'etend en une section de~$\eps_P^*\mathcal L_E$.
Ce sont des m\'etriques alg\'ebriques $v$-adiques; nous noterons
$\bPic_{\alg}(X)$
le groupe des classes d'isom\'etrie de fibr\'es inversibles sur $X_{\bar F}$,
munis de telles m\'etriques.

En particulier, tout diviseur effectif~$D$ sur~$X$ d\'efinit une m\'etrique
sur le fibr\'e en droite $\mathscr O_X(D)$, cette m\'etrique
est d\'efinie par l'adh\'erence de~$D$ dans un mod\`ele $\mathscr X_E$,
o\`u $E$ est une extension de~$F$ de sorte que $D$ soit somme
de points $E$-rationnels de~$X$. Cette d\'efinition s'\'etend
au cas de diviseurs non n\'ecessairement effectifs.
Ce sous-groupe, image de $\Div(X)$ dans~$\bPic_{\alg}(X)$
a un suppl\'ementaire, fourni par les diviseurs verticaux \`a coefficients
rationnels sur un mod\`ele~$\mathcal X_E$.

\`A un $\Q$-diviseur  vertical, on attache une fonction continue,
lin\'eaire par morceaux sur~$R$, dont la  valeur en un point est,
\`a un facteur de normalisation pr\`es, le coefficient de la composante
verticale correspondante. Inversement, une fonction continue~$g$
sur~$R$ d\'efinit une m\'etrique $\norm{\cdot}_g$ sur le fibr\'e en droites
trivial, telle que $\norm{1}_g(x)=\exp(-g(\sp(x)))$. Ce fibr\'e
m\'etris\'e est not\'e $\mathcal O(g)$.

Notons $\bPic(X;R)$ le groupe des (classes d'isom\'etrie) de fibr\'es
m\'etris\'es sur~$X$ qui sont isomorphe \`a un fibr\'e de la forme
$\mathcal O_X(D)\otimes\mathcal O(g)$, o\`u $D$ est un diviseur sur~$X$
et $g$ une fonction continue sur~$\mathcal C(R)$.

Un fibr\'e en droites m\'etris\'e
$\bar L\in \bPic(X;R)$
a une courbure $C_1(\bar L)$, qui est une distribution sur~$R$,
d\'efinie de sorte que pour tout diviseur~$D$ sur~$X$,
$c_1(\mathcal O_X(D))=\mu_D$,
et pour toute fonction continue~$g$ sur~$R$, $c_1(\mathcal O(g))=-\Delta g$,
o\`u $\Delta$ est le laplacien sur le graphe~$R$ (\cf\cite{zhang93}, A.3
pour la d\'efinition). Si la m\'etrique sur $\bar L$ est semi-positive,
cette distribution est une mesure positive; si elle est int\'egrable,
c'est une mesure (\og sign\'ee\fg).
Si $L$ est ample,  on notera
$\mu_{\bar L,R}$ la mesure de masse totale~$1$ sur~$R$
d\'efinie par $c_1(\bar L)/\deg L$.

\begin{theo}\label{thm.graph}
Soit $\bar L$ un fibr\'e en droites ample sur~$X$, muni d'une m\'etrique
ad\'elique int\'egrable.  
Supposons qu'en la place finie~$v$ de~$F$, la m\'etrique de~$\bar L$
appartienne \`a $\bPic(X;R)$.
Alors, pour toute suite $(x_n)$ de points distincts de $X(\bar F)$ telle que
$h_{\bar L}(x_n)\ra h_{\bar L}(x)$, on a la convergence
\[ \mu_{x_n} \ra \mu_{\bar L,R} \]
de mesures sur~$R$.
\end{theo}
L'existence d'une telle suite entra\^{\i}ne en particulier que la mesure
$\mu_{\bar L,R}$ est positive.
\begin{proof}
Soit $\phi$ une fonction lisse sur le graphe~$R$.
Pour tout~$\eps$, le fibr\'e m\'etris\'e
$\bar L(\eps\phi):=\bar L\otimes \mathcal O(\eps\phi)$ est
int\'egrable. 
Il r\'esulte alors de la proposition~\ref{prop.autissier} que
\begin{equation}
\label{eq.zhang2}
 \liminf h_{\bar L(\eps\phi)}(x_n) \geq h_{\bar L(\eps\phi)}(X). 
\end{equation}
Les d\'efinitions des hauteurs entra\^{\i}nent que
\[ h_{\bar L(\eps\phi)}(x) = h_{\bar L}(x) + \int_R \phi \mu_x \log N_v\]
et
\begin{align*}
 h_{\bar L(\eps\phi)}(X) & = \frac{1}{2\deg L} \big(\hc_1(\bar
L(\eps\phi))^2|X\big) \\
& = h_{\bar L}(X) + \eps\log N_v \int_R \phi \frac{c_1(\bar L)}{\deg L}
 + \eps^2 \log N_v\frac{1}{2\deg L} \int_R \phi\Delta\phi. \end{align*}
Par le  m\^eme raisonnement que  dans les autres th\'eor\`emes d'\'equidistribution
de cet article, on a alors
\[ \int _R \phi \mu_{x_n} \ra \int_R \phi \frac{c_1(\bar L)}{\deg L}. \]
Un argument d'approximation imm\'ediat entra\^{\i}ne le r\'esultat pour
toute fonction continue sur~$R$, ce qui cl\^ot la d\'emonstration du th\'eor\`eme.
\end{proof}

L'exemple le plus int\'eressant est peut-\^etre celui d'une
courbe elliptique \`a mauvaise r\'eduction en~$v$. Son graphe
de r\'eduction~$R$ est alors (isom\'etrique \`a) un cercle de
longueur~$\ell=v(\Delta_X)$.
D'un point de vue tr\`es explicite, et \`a la suite
de Tate~\cite{tate95}, une telle courbe admet une uniformisation $p$-adique
comme un quotient rigide analytique de~$\C_p^*$ par un sous-groupe de la
forme $q^\Z$, o\`u $q\in F_v^*$ est de valeur absolue~$<1$. 
L'application de sp\'ecialisation $X(\bar F_v)\ra R$ s'identifie
alors \`a l'application
\[ \log\abs{\cdot}_v \colon \C_p^*\ra\R, \quad\text{induisant}\quad
     \C_p^*/q^\Z \ra \R/(\log\abs {q}_v)\Z = R. \]
La hauteur de N\'eron-Tate sur~$X$ est donn\'ee par le fibr\'e en droites
m\'etris\'e 
\[\bar L= \mathcal O_X(0)\otimes\mathcal O(g_0), \qquad
         g_0(t) = \frac{1}{2\ell} t^2 - \frac12 t + \frac{\ell}{12},
\quad 0\leq t\leq \ell. \]
Ainsi, $c_1(\bar L)=\mathrm dt/\ell$ est la mesure de probabilit\'e
invariante par translation sur le cercle.

\begin{coro}
Supposons que $X$ soit une courbe elliptique \`a mauvaise r\'eduction en~$v$.
Si $(x_n)$ est une suite de points distincts de $X(\bar F)$ 
dont les hauteurs de N\'eron-Tate tendent vers~$0$, les mesures
$\mu_{x_n}$ sur le cercle de r\'eduction convergent vers la mesure
de probabilit\'e invariante.
\end{coro}
En termes imag\'es, les sp\'ecialisations des~$x_n$ balayent
r\'eguli\`erement toutes les composantes irr\'eductibles
de la fibre sp\'eciale de la limite inductive
des mod\`eles de N\'eron.

\subsection{}
Ces r\'esultats sont reli\'es \`a un th\'eor\`eme de Szpiro et Ullmo
dans~\cite{szpiro-u1999}. Si $E$ est une extension finie de~$F$
et $P\in X(E)$, la condition que $P-0+\Phi_P$ ait une intersection nulle
avec tout diviseur vertical d\'efinit un $\Q$-diviseur vertical $\Phi_P$,
unique \`a l'addition d'un l'image inverse d'un diviseur sur~$\Spec\mathfrak
o_E$ pr\`es.
Le nombre d'intersection $-(\Phi_P)^2/[E:\Q]$ est alors positif 
et ne d\'epend pas de l'extension~$E$. 
Szpiro et Ullmo d\'emontrent (\emph{loc.~cit.}, Th\'eor\`eme 1.2)
que pour un point~$P$ d'ordre~$n$, 
\[ -\frac{1}{[F(P):F]}(\Phi_P)^2 = \frac16 \log N_{F/\Q}(\Delta_X) + \mathrm o(d(n)/n^2) . \]  
(Le terme d'erreur est tr\`es souvent nul, \cf\emph{loc.cit.})
Du point de vue des graphes m\'etris\'es, on peut v\'erifier que la contribution
de la place finie~$v$ au membre de gauche, sans son facteur $\log N_v$,
est \'egale \`a
\[ \frac1\ell \int_R t(\ell-t) \,\mathrm dt  = \ell\big(\frac12 -
\frac13\big) = \frac16\ell = \frac16 v(\Delta_X). \]

Plus g\'en\'eralement, le th\'eor\`eme d'\'equidistribution~\ref{thm.graph}
entra\^{\i}ne que pour toute 
suite $(P_n)$ de points distincts de~$X(\bar F)$ dont
les hauteurs de N\'eron-Tate convergent vers~$0$,
$(-\Phi_{P_n})^2/[F(P_n):F]$ converge vers $\frac16 \log
N_{F/\Q}(\Delta_X)$.
En fait, ce dernier point requiert pr\'ecis\'ement une extension  triviale
du th\'eor\`eme~\ref{thm.graph}
dans lequel toutes les places de mauvaise r\'eduction seraient prises en
compte.

\subsection{Minoration de hauteurs}
Comme me l'a sugg\'er\'e P.~Autissier, et ainsi qu'il l'a
lui-m\^eme montr\'e aux places archim\'ediennes, on peut utiliser
ces techniques de changement de m\'etriques pour d\'emontrer
un \'enonc\'e {\og dual\fg}, \`a savoir minorer la limite des hauteurs d'une suite
de points distincts dont les mesures ne convergent pas vers la mesure
attendue.

Un r\'esultat semblable vaut dans le cadre
des espaces de Berkovich, 
avec la m\^eme d\'emonstration que ci-dessous.
Les techniques qu'A.~Thuillier d\'eveloppe dans sa th\`ese
devraient permettre d'obtenir des exemples all\'echants
sur des courbes arbitraires.
Faute d'\^etre capable d'expliciter ces exemples,
je me contente de l'exposer les r\'esultats que l'on peut
obtenir dans le contexte des graphes de r\'eduction.

Soit $C$ une partie ferm\'ee du graphe de r\'eduction~$R$
et soit $(x_n)$ une suite de points distincts de~$X(\bar F)$.
Soit $\bar L$ un fibr\'e inversible sur~$C$ muni d'une m\'etrique ad\'elique
int\'egrable.
On suppose que les mesures $\mu_{x_n}$ convergent vers une mesure~$\nu$
dont le support est contenu dans~$C$; c'est par exemple le cas
si les sp\'ecialisations des conjugu\'es des~$x_n$ appartiennent toutes \`a~$C$.
Soit $\phi$ une fonction continue sur~$R$ telle que $\Delta\phi$
soit une mesure et qui est n\'egative ou nulle sur~$C$.
On a
\[ \liminf h_{\bar L(\phi)}(x_n) = \liminf h_{\bar L}(x_n)
           + \nu(\phi) \log N_v \leq \liminf h_{\bar L}(x_n). \]
Par ailleurs, le th\'eor\`eme de Hilbert-Samuel arithm\'etique entra\^{\i}ne
(\emph{bis repetita\dots}) l'in\'egalit\'e
\[ \liminf h_{\bar L(\phi)}(x_n) \geq h_{\bar L(\phi)}(X)
   =  h_{\bar L}(X) 
+ \left(\int_R \phi \frac{c_1(\bar L)}{\deg L} + \frac1{2\deg L} \int_R
\phi\Delta\phi\right)\log N_v. \]
On a donc une minoration
\begin{equation}
\label{eq.minorh}
 \liminf h_{\bar L}(x_n) \geq h_{\bar L}(X) + 
\left( \int_R \phi \frac{c_1(\bar
L)}{\deg L} - \frac1{2\deg L} \int_R (\phi')^2\right)
\log N_v. \end{equation}

L'optimisation de cette in\'egalit\'e conduit, si c'est possible,
\`a prendre pour~$\phi$ une fonction nulle sur~$C$ et telle que $\Delta\phi$
soit \'egale \`a $c_1(\bar L)/\deg L$ hors de~$C$.

Contentons nous d'expliciter ici  le cas o\`u $X$ est une courbe
elliptique \`a mauvaise r\'eduction, le graphe~$R$ est donc un cercle
de longueur~$\ell$, et supposons que le compl\'ementaire de~$C$
dans~$R$ soit une r\'eunion finie d'intervalles disjoints
$\mathopen]a_i,b_i\mathclose[$, pour $1\leq i\leq t$,
avec $0\leq a_1<b_1\leq \dots\leq a_t<b_t\leq \ell$.
Sur $\mathopen]a_i,b_i\mathclose[$, on pose
$\phi(t)=c_i (t-a_i)(b_i-t)$, o\`u $c_i$
est une constante \`a d\'eterminer; sur~$C$, on pose $\phi(t)=0$.
Sur $\mathopen]a_i,b_i\mathclose[$, on a ainsi $\Delta\phi=-2c_i$,
si bien que
le membre de droite de l'\'equation~\eqref{eq.minorh}
vaut
\[
\frac{\log N_v}\ell \sum_{i=1}^t \int_{a_i}^{b_i} c_i(1-\ell c_i)(t-a_i)(b_i-t)
= \frac{\log N_v}{6\ell} \sum_{i=1}^t (b_i-a_i)^3 c_i(1-\ell c_i) . \] 
On choisit $c_i=1/2\ell$, d'o\`u la minoration
\[ \liminf h_{\bar L}(x_n) \geq \frac1{24\ell^2} \sum_{i=1}^t (b_i-a_i)^3
\log N_v.
\]

Soyons encore plus explicites en donnant trois exemples.

1) Si l'on demande aux points~$x_n$ de passera par la composante neutre
du mod\`ele de N\'eron en~$v$, le support de la mesure~$\nu$ est le point du 
graphe~$R$ correspondant \`a cette composante:  on a donc $t=1$,
$a_1=0$ et $b_1=\ell$. L'in\'egalit\'e pr\'ec\'edente devient
\[ \liminf h_{\bar L}(x_n) \geq \frac1{24}  v(\Delta_X) \log N_v. \]
Modulo l'extension  de ce th\'eor\`eme
\`a un nombre fini de places, on voit
que pour toute suite $(x_n)$ de points distincts de~$X(\bar F)$
qui sont des points entiers de la composante neutre du mod\`ele
de N\'eron de~$X$ sur~$\mathfrak o_F$, la liminf des hauteurs
de N\'eron-Tate des~$x_n$ est au moins $\log\Delta_X / 24$.
On retrouve ainsi un th\'eor\`eme de Ullmo (\cite{ullmo1995}, th\'eor\`eme 4.1),
\`a un facteur~$2$ pr\`es qui provient de la normalisation de la hauteur
de N\'eron-Tate dans~\loccit.

2) Si l'on demande aux points~$x_n$ d'\^etre des points entiers
du mod\`ele de N\'eron (suppos\'e semi-stable) de~$X$ sur le corps~$F$,
cela revient \`a prendre $t=\ell$, et $a_i=i-1$, $b_i=i$
pour $1\leq i\leq t$. Alors,
\[ \liminf h_{\bar L}(x_n) \geq \frac1{24 v(\Delta_X) }\log N_v. \]

3) On peut aussi exiger seulement que les points~$x_n$ \'evitent un point
singulier de la fibre en~$v$ du mod\`ele minimal r\'egulier de~$X$ sur~$F$.
Cela revient \`a prendre $t=1$, $a_1=0$ et $b_1=1$.
On obtient alors
\[ \liminf h_{\bar L}(x_n) \geq \frac1{24 v(\Delta_X)^2}\log N_v. \]

Notons que ces exemples montrent que le th\'eor\`eme 2.2 d'Autissier
dans~\cite{autissier2001} n'est pas valable sans une hypoth\`ese
de bonne r\'eduction.

\bibliographystyle{smfplain}
\bibliography{aclab,acl,equi}

\providecommand{\noopsort}[1]{}\providecommand{\url}[1]{\textit{#1}}
\providecommand{\bysame}{\leavevmode ---\ }
\providecommand{\og}{``}
\providecommand{\fg}{''}
\providecommand{\smfandname}{\&}
\providecommand{\smfedsname}{\'eds.}
\providecommand{\smfedname}{\'ed.}
\providecommand{\smfmastersthesisname}{M\'emoire}
\providecommand{\smfphdthesisname}{Th\`ese}
\begin{thebibliography}{10}

\bibitem{abbes-b95}
{\scshape A.~Abbes {\normalfont \smfandname} T.~Bouche} -- {\og
  Th{\'e}or{\`e}me de {H}ilbert--{S}amuel {\og arithm{\'e}tique\fg}\fg},
  \emph{Ann. Inst. Fourier (Grenoble)} \textbf{45} (1995), no.~2, p.~375--401.

\bibitem{autissier2001}
{\scshape P.~Autissier} -- {\og Points entiers et th\'eor\`emes de {B}ertini
  arithm\'etiques\fg}, \emph{Ann. Inst. Fourier (Grenoble)} \textbf{51} (2001),
  no.~6, p.~1507--1523, \emph{corrigendum:} \textbf{52} (2002), no. 1,
  p.~303--304.

\bibitem{autissier2001b}
\bysame , {\og Points entiers sur les surfaces arithm\'etiques\fg}, \emph{J.
  Reine Angew. Math.} \textbf{531} (2001), p.~201--235.

\bibitem{baker-h2003}
{\scshape M.~Baker {\normalfont \smfandname} L.-C. Hsia} -- {\og Canonical
  heights, transfinite diameters, and polynomial dynamics\fg},
  arXiv:math.NT/0305181.

\bibitem{baker-r2004}
{\scshape M.~Baker {\normalfont \smfandname} R.~Rumely} -- {\og
  {Equidistribution of small points, rational dynamics, and potential
  theory}\fg}, arXiv, math.NT/0407426.

\bibitem{berkovich1990}
{\scshape V.~G. Berkovich} -- \emph{Spectral theory and analytic geometry over
  non-{A}rchimedean fields}, Mathematical Surveys and Monographs, vol.~33,
  American Mathematical Society, Providence, RI, 1990.

\bibitem{bilu97}
{\scshape {\relax Yu}.~Bilu} -- {\og Limit distribution of small points on
  algebraic tori\fg}, \emph{Duke Math. J.} \textbf{89} (1997), no.~3,
  p.~465--476.

\bibitem{bosch-l1993}
{\scshape S.~Bosch {\normalfont \smfandname} W.~L{\"u}tkebohmert} -- {\og
  Formal and rigid geometry. {I}. {R}igid spaces\fg}, \emph{Math. Ann.}
  \textbf{295} (1993), no.~2, p.~291--317.

\bibitem{bost-g-s94}
{\scshape J.-B. Bost, H.~Gillet {\normalfont \smfandname} C.~Soul{\'e}} -- {\og
  Heights of projective varieties and positive {G}reen forms\fg}, \emph{J.
  Amer. Math. Soc.} \textbf{7} (1994), p.~903--1027.

\bibitem{david-p98}
{\scshape S.~David {\normalfont \smfandname} P.~Philippon} -- {\og Minorations
  des hauteurs normalis{\'e}es des sous-vari{\'e}t{\'e}s de vari{\'e}t{\'e}s
  ab{\'e}liennes\fg}, in \emph{International Conference on Discrete Mathematics
  and Number Theory} (Tiruchirapelli, 1996), Contemp. Math., 1998, p.~333--364.

\bibitem{faltings84}
{\scshape G.~Faltings} -- {\og Calculus on arithmetic surfaces\fg}, \emph{Ann.
  of Math.} \textbf{119} (1984), p.~387--424.

\bibitem{favre-r-t2004}
{\scshape C.~Favre {\normalfont \smfandname} J.~Rivera-Letelier} -- {\og
  {Equidistribution des points de petite hauteur}\fg}, arXiv, math.NT/0407471.

\bibitem{gillet-s88}
{\scshape H.~Gillet {\normalfont \smfandname} C.~Soul{\'e}} -- {\og Amplitude
  arithm{\'e}tique\fg}, \emph{C. R. Acad. Sci. Paris S{\'e}r. I Math.}
  \textbf{307} (1988), p.~887--890.

\bibitem{gillet-s90}
\bysame , {\og Arithmetic intersection theory\fg}, \emph{Publ. Math. Inst.
  Hautes {\'E}tudes Sci.} \textbf{72} (1990), p.~94--174.

\bibitem{gillet-s92}
\bysame , {\og An arithmetic {R}iemann--{R}och theorem\fg}, \emph{Invent.
  Math.} \textbf{110} (1992), p.~473--543.

\bibitem{gubler1997}
{\scshape W.~Gubler} -- {\og Heights of subvarieties over {$M$}-fields\fg}, in
  \emph{Arithmetic geometry} (F.~Catanese, \smfedname), Symp. Math., vol.~37,
  1997, p.~190--227.

\bibitem{gubler1998}
\bysame , {\og Local heights of subvarieties over non-archimedean fields\fg},
  \emph{J. Reine Angew. Math.} \textbf{498} (1998), p.~61--113.

\bibitem{gubler2003}
\bysame , {\og Local and canonical heights of subvarieties\fg}, \emph{Ann.
  Scuola Norm. Sup. Pisa} \textbf{2} (2003), no.~4, p.~711--760.

\bibitem{maillot2000}
{\scshape V.~Maillot} -- {\og G\'eom\'etrie d'{A}rakelov des vari\'et\'es
  toriques et fibr\'es en droites int\'egrables\fg}, \emph{M{\'e}m. Soc. Math.
  France} (2000), no.~80, p.~129.

\bibitem{mainetti2001}
{\scshape N.~Ma{\"{\i}}netti} -- {\og Metrizability of some analytic affine
  spaces\fg}, in \emph{$p$-adic functional analysis (Ioannina, 2000)}, Lecture
  Notes in Pure and Appl. Math., vol. 222, Dekker, New York, 2001, p.~219--225.

\bibitem{rumely1989}
{\scshape R.~Rumely} -- \emph{Capacity theory on algebraic curves}, Lecture
  Notes in Math., vol. 1378, Springer-Verlag, Berlin, 1989.

\bibitem{rumely1995}
\bysame , {\og An intersection pairing for curves, with analytic contributions
  from non-{A}rchimedean places\fg}, in \emph{Number theory (Halifax, NS,
  1994)}, CMS Conf. Proc., vol.~15, Amer. Math. Soc., Providence, RI, 1995,
  p.~325--357.

\bibitem{szpiro-u1999}
{\scshape L.~Szpiro {\normalfont \smfandname} E.~Ullmo} -- {\og Variation de la
  hauteur de {F}altings dans une classe de {$\overline{\mathbf Q}$}-isog\'enie
  de courbe elliptique\fg}, \emph{Duke Math. J.} \textbf{97} (1999), no.~1,
  p.~81--97.

\bibitem{szpiro-u-z97}
{\scshape L.~Szpiro, E.~Ullmo {\normalfont \smfandname} S.-W. Zhang} -- {\og
  {\'E}quidistribution des petits points\fg}, \emph{Invent. Math.} \textbf{127}
  (1997), p.~337--348.

\bibitem{tate95}
{\scshape J.~Tate} -- {\og A review of non-archimedean elliptic function\fg},
  in \emph{Conference on elliptic curves and modular forms} (Hong Kong, 1993),
  1995, p.~162--184.

\bibitem{ullmo1995}
{\scshape E.~Ullmo} -- {\og Points entiers, points de torsion et amplitude
  arithm\'etique\fg}, \emph{Amer. J. Math.} \textbf{117} (1995), no.~4,
  p.~1039--1055.

\bibitem{zhang92}
{\scshape S.-W. Zhang} -- {\og Positive line bundles on arithmetic
  surfaces\fg}, \emph{Ann. of Math.} \textbf{136} (1992), no.~3, p.~569--587.

\bibitem{zhang93}
\bysame , {\og Admissible pairing on a curve\fg}, \emph{Invent. Math.}
  \textbf{112} (1993), no.~1, p.~171--193.

\bibitem{zhang95}
\bysame , {\og Positive line bundles on arithmetic varieties\fg}, \emph{J.
  Amer. Math. Soc.} \textbf{8} (1995), p.~187--221.

\bibitem{zhang95b}
\bysame , {\og Small points and adelic metrics\fg}, \emph{J. Algebraic
  Geometry} \textbf{4} (1995), p.~281--300.

\bibitem{zhang98}
\bysame , {\og Equidistribution of small points on abelian varieties\fg},
  \emph{Ann. of Math.} \textbf{147} (1998), no.~1, p.~159--165.

\end{thebibliography}
\end{document}